\documentclass[11pt]{amsart}

\usepackage{amsmath}
\usepackage{amssymb}
\usepackage{sty1}

\makeatletter
\@addtoreset{equation}{section}
\makeatother

\title[Representations of Group Algebras]
{Representations of Group Algebras in Spaces of Completely
Bounded Maps}

\author{Roger R.\ Smith and Nico Spronk}

\def\augidealg{\mathrm{I}_0(G)}
\def\augideal#1{\mathrm{I}_0(#1)}
\def\bdop#1{\mathcal{B}(#1)}
\def\bdophil{\mathcal{B}(\mathcal{H})}
\def\bfalp{\pmb{\alpha}}
\def\bflam{\pmb{\lambda}}
\def\bfeta{\pmb{\eta}}
\def\bfxi{\pmb{\xi}}
\def\blinfty#1{\mathrm{L}^\infty(#1)}
\def\blone#1{\mathrm{L}^1(#1)}
\def\bloneg{\mathrm{L}^1(G)}

\def\cbop#1{\mathcal{CB}(#1)}

\def\ceeo#1{c_{\,0}(#1)}
\def\ceeos{c_{\,0}\text{-}}
\def\coef#1{\mathrm{F_{#1}}}
\def\contb#1{\mathcal{C}_b(#1)}
\def\conto#1{\mathcal{C}_0(#1)}
\def\cstarg{\mathrm{C}^*(G)}
\def\cstardpi{\mathrm{D}^*_{\pi}}
\def\cstarpi{\mathrm{C}^*_\pi}
\def\cstarpih{\mathrm{C}^*_{\pi|_H}}
\def\dgrp{\what{G}}
\def\dgrph{\what{H}}
\def\gampi{\Gamma_\pi}
\def\gampih{\Gamma_{\pi|_H}}

\def\hil{\mathcal{H}}

\def\tiljrpi{\til{\jmath}_{r_\pi}}
\def\loneg{\ell^1(G)}
\def\ltwo#1{\ell^2(#1)}
\def\matn#1{\mathrm{M}_n(#1)}
\def\meas#1{\mathrm{M}(#1)}
\def\measg{\mathrm{M}(G)}
\def\meash{\mathrm{M}(H)}
\def\mstarpi{\mathrm{M}^*_\pi}
\def\mstarpih{\mathrm{M}^*_{\pi|_H}}
\def\nlift#1{{#1}^{(n)}}

\def\supp#1{\mathrm{supp}(#1)}
\def\tilfpi{\til{F}_{\pi}}
\def\tilfpih{\til{F}_{\pi|_H}}

\def\transpose#1{{#1}^{\mathrm{t}}}
\def\unop#1{\mathcal{U}(#1)}
\def\unophil{\mathcal{U}(\mathcal{H})}
\def\val#1{\mathrm{V}^\infty(#1)}
\def\valo#1{\mathrm{V}^0(#1)}
\def\valb#1{\mathrm{V}^b(#1)}
\def\valoc#1{\mathrm{V}_0(#1)}
\def\vnpi{\mathrm{VN}_\pi}
\def\wbdop#1{\mathcal{B}^\sigma(#1)}
\def\wcbop#1{\mathcal{CB}^\sigma(#1)}

\def\tens{\otimes}
\def\ehnorm#1{\left\|#1\right\|_{eh}}
\def\hnorm#1{\left\|#1\right\|_h}
\def\htens{\otimes^h}
\def\ehtens{\otimes^{eh}}
\def\whtens{\otimes^{w^*h}}
\def\whnorm#1{\left\|#1\right\|_{w^*h}}

\def\measalp{\mathrm{M}_\alp}
\def\cmeasalp{\mathrm{C}_\alp}
\def\dmeasalp{\mathrm{D}_\alp}
\def\inv{{\mathrm{inv}}}

\begin{document}
\maketitle

\newtheorem{ehtpisalg}{Proposition}[section]
\newtheorem{posiscp}[ehtpisalg]{Lemma}
\newtheorem{cpmapchar}[ehtpisalg]{Lemma}
\newtheorem{posabelian}[ehtpisalg]{Theorem}

\newtheorem{alpishomo}{Proposition}[section]
\newtheorem{moduleexample}[alpishomo]{Example}
\newtheorem{repwunit}[alpishomo]{Proposition}
\newtheorem{repwunit1}[alpishomo]{Corollary}
\newtheorem{repsingeneral}[alpishomo]{Proposition}
\newtheorem{contincb}[alpishomo]{Remark}
\newtheorem{repinhtp}[alpishomo]{Theorem}
\newtheorem{kergampi}[alpishomo]{Proposition}

\newtheorem{abelianreps0}{Lemma}[section]
\newtheorem{abelianreps}[abelianreps0]{Theorem}
\newtheorem{abelianreps1}[abelianreps0]{Corollary}
\newtheorem{abelianreps2}[abelianreps0]{Corollary}
\newtheorem{squareex}[abelianreps0]{Example}
\newtheorem{injectiveex}[abelianreps0]{Question}

\footnote{{\it Date}: \today.


The second author was supported by an NSERC PDF.}

\begin{abstract}
Let $G$ be a locally compact group, $\pi:G\to\unophil$ be a strongly
continuous unitary representation, and $\wcbop{\bdophil}$ the space
of normal completely bounded maps on $\bdophil$.
We study the range of the map
\[
\gampi:\measg\to\wcbop{\bdophil},\quad
\gampi(\mu)= \int_G \pi(s)\tens\pi(s)^*d\mu(s)
\]
where we identify $\wcbop{\bdophil}$ with the extended Haagerup
tensor product $\bdophil\ehtens\bdophil$.  
We use the fact that the C*-algebra generated by integrating $\pi$
to $\bloneg$ is unital exactly when $\pi$ is norm continuous, to show that
$\gampi(\bloneg)\subset\bdophil\htens\bdophil$ exactly when $\pi$
is norm continuous.  For the case that $G$ is abelian, we study
$\gampi(\measg)$ as a subset of the Varopoulos algebra.
We also characterise positive definite elements of the Varopoulos
algebra in terms of completely positive operators.
\end{abstract}

\pagebreak

\section{Introduction}

In \cite{stormer}, St{\o}rmer conducted an interesting study of spaces
of completely bounded maps on $\bdophil$.  For subalgebras $\fA$ and $\fB$
of $\bdophil$ he defined what is now known as the Haagerup tensor
product $\fA\htens\fB$, as a completion of the set of elementary operators
of the form $x\mapsto\sum_{i=1}^n a_ixb_i$ where each $a_i\in\fA$ and
each $b_i\in\fB$.  This approach gives the same tensor product norm
as that in the more standard approach (see \cite{effrosrB}, for example),
as shown in \cite{smith}.

If $G$ is an abelian group and $\pi:G\to\unophil$ is a strongly continuous
unitary representation, the homomorphism $\gampi$
from the measure algebra
$\measg$ to the space $\wcbop{\bdophil}$ of
normal completely bounded maps on $\bdophil$, defined by
\begin{equation}\label{eq:themap}
\gampi(\mu)=\int_G\pi(s)\tens\pi(s)^*d\mu(s)
\end{equation}
was studied by St{\o}rmer.  (We identify $\wcbop{\bdophil}$ with
the extended Haagerup tensor product $\bdophil\ehtens\bdophil$
from \cite{blechers} and \cite{effrosr}.)
He used this homomorphism to generate many examples
of regular and non-regular Banach subalgebras of $\wcbop{\bdophil}$.
It was shown in  \cite[Lem. 5.6]{stormer} that if 
$\pi$ is norm continuous (i.e.\ continuous when the norm topology is 
placed on $\unophil$) then for any $f\iin\bloneg$
\begin{equation}\label{eq:uhoh}
\gampi(f)=\int_G f(s)\pi(s)\tens\pi(s)^*ds\in\cstarpi\htens\cstarpi
\end{equation}
where $\cstarpi$ is the C*-algebra generated by
$\cbrac{\int_G f(s)\pi(s)ds:f\in\bloneg}$.  

We note that for an arbitrary locally compact group $G$, the map
$\Gamma_\lam$ as in (\ref{eq:themap}), where $\lam$ is the left regular 
representation, was studied in \cite{ghahramani} and \cite{neufang}.

\medskip
In this paper we will make use of the theory of completely bounded normal
maps on $\bdophil$ from \cite{smith} to study the range of $\gampi$.
We show that, for a general locally compact group G,
\[
\gampi(\bloneg)\subset\cstarpi\ehtens\cstarpi
\]
where $\ehtens$ denotes the extended Haagerup tensor product from
\cite{effrosr}, \cite{effroskr} and \cite{blechers}.
Moreover, using the fact that $\cstarpi$ is unital exactly when
the representation $\pi:G\to\unophil$ is norm continuous, we show
that the validity of (\ref{eq:uhoh}) for every $f$ in $\bloneg$
gives a characterisation of the norm continuity of $\pi$.

In the case that $G$ is abelian, we develop the
``Fourier-Stieltjes transform'' for $\gampi(\measg)$.  The range of this
transform is a Varopoulos type algebra $\valb{E_\pi}$, which will be 
defined below.  We use some general results on completely positive maps
to characterise complete positivity of elements of $\valb{E_\pi}$,
as operators on $\bdophil$, extending some results from \cite{stormer}. 
In particular, we characterise 
those $\mu\iin\measg$ for which $\gampi(\mu)$ is completely positive.

\section{Spaces of Normal Completely Bounded Maps}

Let $\hil$ be a Hilbert space, let $\bdophil$ be the space of bounded 
operators 
on $\hil$ and let $\fV$ and $\fW$ be closed subspaces of $\bdophil$.
The {\it Haagerup tensor product} $\fV\htens\fW$ is defined in 
\cite{haagerup} and \cite{effrosk}.  The {\it extended Haagerup tensor product}
$\fV\ehtens\fW$
is developed in \cite{effrosr} and \cite{effroskr}; and also in 
\cite{blechers}, but in the context of dual spaces where it is called the
``weak* Haagerup tensor product'' and denoted $\fV\whtens\fW$.  It is shown
in \cite{spronk} that the approach of \cite{blechers} can be modified
to develop the extended Haagerup tensor product in general.

Following \cite{spronk}, we thus define $\fV\ehtens\fW$ to be the space
of all (formal) series $\sum_{i\in I}v_i\tens w_i$ where each $v_i\in\fV$, each
$w_i\in\fW$, and each of the series $\sum_{i\in I}v_iv_i^*$ and 
$\sum_{i\in I}w_i^*w_i$
converges weak* in $\bdophil$.  The index set $I$ is established to have
cardinality $|I|=\dim\hil$.  Two series $\sum_{i\in I}v_i\tens w_i$
and $\sum_{i\in I}v_i'\tens w_i'$ define the same element of $\fV\ehtens\fW$
provided $\sum_{i\in I}v_ix w_i=\sum_{i\in I}v_i'x w_i'$ for each 
$x\iin\bdophil$.  Then $\fV\ehtens\fW$ is a Banach space when endowed with 
the norm
\[
\ehnorm{T}=\inf\cbrac{
\norm{\sum_{i\in I}v_iv_i^*}^{1/2}\norm{\sum_{i\in I}w_i^*w_i}^{1/2}:
T=\sum_{i\in I}v_i\tens w_i}
\]
and the infimum is attained.  As in \cite{blechers},
note that the Haagerup tensor product 
$\fV\htens\fW$ may be realized as the set of those $T$ in $\fV\ehtens\fW$
which admit a representation $T=\sum_{i\in I}v_i\tens w_i$ where
$\sum_{i\in I}v_iv_i^*$ and $\sum_{i\in I}w_i^*w_i$ converge in norm.
It is easy to see that any element $T$ of $\fV\htens\fW$ may thus be written
with a countable index set
as $T=\sum_{i=1}^\infty v_i\tens w_i$.

The space $\fV\ehtens\fW$ has two natural, though more extrinsic
descriptions.
First, if $\fV$ and $\fW$ are each weak* closed subspaces of $\bdophil$,
they have respective preduals $\fV_*$ and $\fW_*$.  For example,
\[
\fV_*=\bdophil_*/\{\omega\in\bdophil_*:\omega(v)=0\text{ for all }
v\iin\fV\} 
\]
which is an operator space when endowed with the quotient 
structure from the predual operator space structure on $\bdophil_*$.
Then $\fV\ehtens\fW$ is the dual space of $\fV_*\htens\fW_*$ via the pairing
\begin{equation}\label{eq:htpdual}
\dpair{\sum_{i\in I}v_i\tens w_i}{\omega\tens\nu}
=\sum_{i\in I}\omega(v_i)\nu(w_i).
\end{equation}
A proof of this can be found in \cite{blechers} or \cite{effrosr}.
In particular, 
$\bdophil\ehtens\bdophil\cong\brac{\bdophil_*\htens\bdophil_*}^*$.

Let $\wcbop{\bdophil}$ denote the space of normal completely bounded operators
on $\bdophil$.
The map $\theta:\bdophil\ehtens\bdophil\to\wcbop{\bdophil}$ given by
\[
\theta\brac{\sum_{i\in I}v_i\tens w_i}x=\sum_{i\in I}v_ix w_i,
\ffor x\iin\bdop{\hil}
\]
is a surjective isometry by \cite{haagerup} or \cite{smith}.  Moreover,
$\theta$ is still an isometry when restricted to the spaces
$\fV\ehtens\fW$ or $\fV\htens\fW$.  For notational ease we will 
simply identify $\fV\ehtens\fW$ and $\fV\htens\fW$ as subspaces of
$\wcbop{\bdophil}$ in the sequel, and omit the map $\theta$.
In particular, we view $\bdophil\htens\bdophil$ as being the completion 
in the completely bounded operator norm of
the space of elementary operators $x\mapsto\sum_{i=1}^nv_ixw_i$
on $\bdophil$.  The composition of operators in
$\wcbop{\bdophil}$ induces a product in $\bdophil\ehtens\bdophil$, making it 
a Banach algebra.  This product is given on elementary tensors by
\[
(a\tens b)\comp(c\tens d)=ac\tens db.
\]
The following is an extension of a theorem from \cite{blecher}, whose
proof is much like the one offered there.

\begin{ehtpisalg}\label{prop:ehtpisalg}
If $\fA$ and $\fB$ are norm closed subalgebras of $\bdophil$, then
$\fA\ehtens\fB$ is a subalgebra of $\bdophil\ehtens\bdophil$.
If $\fV$ is a (left) $\fA$-module and $\fW$ is a (right) $\fB$-module
in $\bdophil$, then $\fV\ehtens\fW$ is a (left) $\fA\ehtens\fB$-module 
in $\bdophil\ehtens\bdophil$.
\end{ehtpisalg}


If $\Omega\in\bdophil^*$ then the {\it left} and {\it right slice maps}
$L_\Omega,R_\Omega:\bdophil\ehtens\bdophil\to\bdophil$ are given for
$T=\sum_{i\in I}v_i\tens w_i$ by
\begin{equation}\label{eq:slicemaps}
L_\Omega T=\sum_{i\in I}\Omega(v_i)w_i\quad\aand\quad
R_\Omega T=\sum_{i\in I}\Omega(w_i)v_i.
\end{equation}
These series each converge in norm as is shown in 
\cite[Thm.\ 2.2]{spronk}.  Moreover, it is shown there
that for any pair of closed subspaces $\fV$ and $\fW$ 
of $\bdophil$, $\fV\ehtens\fW$ consists exactly of those $T$ in 
$\bdophil\ehtens\bdophil$ for which $L_\Omega T\in\fW$ and $R_\Omega T\in\fV$
for each $\Omega\iin\bdophil^*$ (or for which $L_\ome T\in\fW$ and 
$R_\ome T\in\fV$ for each $\ome\iin\bdophil_*$).  These results extend 
\cite[Thm. 4.5]{smith}.

We will finish this section with a theorem on completely positive maps which
will be useful in Section \ref{sec:abelian}.  We will first need
some general preliminary results which are modeled on 
results from \cite{smith}.

A closed subalgebra $\fB$ of $\bdophil$ is called {\it locally cyclic} 
if for each finite dimensional subspace $\fL$ of $\hil$, there
is a vector $\xi\iin\hil$ such that $\wbar{\fB\xi}\supset\fL$.
We note, for example, that if $\fB$ is a maximal abelian self-adjoint
subalgebra of $\bdophil$ then it is locally cyclic.  Indeed if 
$\xi_1,\dots,\xi_n$ span
$\fL$, consider the orthogonal projections $p_1,p_2,\dots,p_n$
whose respective ranges are
\[
\wbar{\fB\xi_1},\;
\wbar{\fB\xi_2}\ominus\wbar{\fB\xi_1},\;
\dots,\;\wbar{\fB\xi_n}\ominus\bigoplus_{i=1}^{n-1}\wbar{\fB\xi_i}.
\]
Then each $p_i\in\fB'=\fB$, and
$\xi=\xi_1+p_2\xi_2+\dots+p_n\xi_n$ satisfies $\wbar{\fB\xi}\supset\fL$.

The following is an adaptation of \cite[Thm.\ 2.1]{smith}.

\begin{posiscp}\label{lem:posiscp}
If $\fB$ is a locally cyclic C*-subalgebra of $\bdophil$
and $T:\bdophil\to\bdophil$ is a positive map which is also a
$\fB$-bimodule map, then $T$ is completely positive.
\end{posiscp}

\proof Let us fix $n$, a positive matrix $[x_{ij}]$ in $\matn{\bdophil}$ 
and a column vector
$\bfxi=\transpose{[\xi_1\cdots\xi_n]}\iin\hil^n$ with $\norm{\bfxi}<1$.
Then, given $\eps>0$, there is vector $\xi$ in $\hil$ and elements
$b_1,\dots,b_n$ in $\fB$ such that the vector $\bfeta
=\transpose{[b_1\xi\cdots b_n\xi]}$ satisfies $\norm{\bfxi-\bfeta}<\eps$ and
$\norm{\bfeta}<1$.  Leting $\nlift{T}:\matn{\bdophil}\to\matn{\bdophil}$
be the amplification of $T$, we have
\begin{align*}
\inprod{\nlift{T}[x_{ij}]\bfeta}{\bfeta}
&=\inprod{[Tx_{ij}]\begin{bmatrix} b_1\xi \\ \vdots \\ b_n\xi \end{bmatrix}}
{\begin{bmatrix} b_1\xi \\ \vdots \\ b_n\xi \end{bmatrix}} \\
&=\sum_{i,j=1}^n \inprod{b_i^*T(x_{ij})b_j\xi}{\xi}
=\inprod{T\brac{\sum_{i,j=1}^nb_i^*x_{ij}b_j}\xi}{\xi}\geq 0
\end{align*}
and
\[
\abs{\inprod{\nlift{T}[x_{ij}]\bfeta}{\bfeta}-
\inprod{\nlift{T}[x_{ij}]\bfxi}{\bfxi}}<\brac{\norm{\nlift{T}}+1}\eps.
\]
Since $\eps$ can be chosen arbitrarily small, we conclude that
$\inprod{\nlift{T}[x_{ij}]\bfxi}{\bfxi}\geq 0$.  Hence $T$ is
completely positive. \endpf

If a family of operators $\{b_i\}_{i\in I}$ from $\bdophil$ defines a 
bounded row matrix $B=[\cdots b_i\cdots]$, i.e.\ $\sum_{i\in I}b_ib_i^*$
converges weak* in $\bdophil$, then the product 
$B\mult\bflam=\sum_{I\in I}\lam_ib_i$ converges in norm and thus defines
an element of $\bdophil$ for each $\bflam=\transpose{[\cdots \lam_i\cdots]}
\iin\ltwo{I}$.  We say that the set $\{b_i\}_{i\in I}$ is {\it strongly
independent} if $B\mult\bflam=0$ only when $\bflam=0$.  This is an obvious
extension of the usual notion of linear independence, and can
be easily adapted to column matrices.  Elements of $\bdophil\ehtens\bdophil$
admit many different representations, and strong independence was introduced 
in \cite{smith} to handle the difficulties caused by this.

The following is an adaptation of \cite[Thm.\ 3.1]{smith}.

\begin{cpmapchar}\label{lem:cpmapchar}
If $\fA$ is a C*-subalgebra of $\bdophil$ and $T\in\fA\ehtens\fA$,
then $T$ is completely positive if and only if there is a 
strongly independent family
$\{a_i\}_{i\in I}$ from $\fA$ for which $\sum_{i\in I}a_ia_i^*$ converges
weak* in $\bdophil$ and $T=\sum_{i\in I}a_i\tens a_i^*$.
\end{cpmapchar}

\proof  We need only to prove that the first condition implies the second.  

If $T$ is completely positive and normal on $\bdophil$, then
its restriction to the algebra of compact operators 
$T|_{\fK(\hil)}$ is a completely
positve map which determines $T$.  Using Stinespring's theorem and
the representation theory for $\fK(\hil)$, just as in 
\cite[Thm.\ 3.1]{smith} or \cite{haagerup}, 
we obtain a family $\{b_j\}_{j\in J}$ from $\bdophil$
for which $\sum_{j\in J}b_jb_j^*$ converges weak* in $\bdophil$
and $T=\sum_{j\in J}b_j\tens b_j^*$.  We see that $J$ can be any
index set whose cardinality coincides with the Hilbertian dimension of $\hil$.
Let $B=[\cdots b_j\cdots]$.

Now we let 
\[
\fL=\{\bflam\in\ltwo{J}:B\mult\bflam=0\} 
\]
and partition
$J=I'\cup I$ in such a way that there is an orthonormal basis
$\{\bflam_j\}_{j\in J}$ of $\ltwo{J}$ for which 
\[
\wbar{\spn}\{\bflam_i\}_{i\in I'}=\fL\quad\aand\quad
\wbar{\spn}\{\bflam_i\}_{i\in I}=\fL^\perp.
\]
Let $U$ denote the $J\cross J$ unitary matrix whose columns are the vectors
$\{\bflam_j\}_{j\in J}$.  Let $A=[\cdots a_j\cdots]=B\mult U$.  Note that
$a_j=0$ for each $j\iin I'$.  Then for any $x\iin\bdophil$, letting
$x^J$ denote the $J\cross J$ diagonal matrix which is the amplification of
$x$, we have
\[
Tx=\sum_{j\in J}b_jxb_j^*=Bx^JB^*=B\mult Ux^JU^*\mult B^*=Ax^JA^*
=\sum_{i \in I}a_ixa_i^*.
\]
We have that $\{a_i\}_{i\in I}$ is strongly independent, for if $\bfalp
=\transpose{[\cdots \alp_i\cdots]}$ in $\ltwo{I}$ is such that
$A\mult\bfalp=0$, then
\[
0=A\mult\bfalp=\sum_{i\in I}\alp_i a_i=\sum_{i\in I}\alp_i B\mult\bflam_i
=B\mult\brac{\sum_{i\in I}\alp_i\bflam_i}
\]
so $\sum_{i\in I}\alp_i\bflam_i\in\fL\cap\fL^\perp$, whence $\bfalp=0$.
Hence
\[
T=\sum_{i \in I}a_i\tens a_i^*
\]
where $\{a_i\}_{i\in I}$ is strongly independent.
It remains to show that $\{a_i\}_{i\in I}\subset\fA$.

Since $\{a_i\}_{i\in I}$ is strongly independent, so too is
$\{a_i^*\}_{i\in I}$.  Hence by \cite[Lem.\ 2.2]{allenss}, the space
\[
\cbrac{\transpose{[\cdots \Omega(a_i^*)\cdots]}:\Omega\in\bdophil^*}
\]
is dense in $\ltwo{I}$.  Thus, given a fixed index $i_0\iin I$, there
is a (not necessarily bounded)
sequence $\{\Omega_n\}_{n=1}^\infty$ from $\bdophil^*$ such that
\[
a_{i_0}=\lim_{n\to\infty}\sum_{i\in I}\Omega_n(a_i^*)a_i=
\lim_{n\to\infty}R_{\Omega_n}T.
\]
Since $R_\Omega T\in\fA$ for each right slice map $R_\Omega$, it follows that
$a_{i_0}\in\fA$. \endpf

If $E$ is any locally compact space we let
\begin{align}\label{eq:valgebras}
\valoc{E}&=\conto{E}\htens\conto{E} \notag \\
\valo{E}&=\conto{E}\ehtens\conto{E} \\
\aand\valb{E}&=\contb{E}\ehtens\contb{E}. \notag
\end{align}
These spaces are discussed in \cite{spronk}.  These all may be regarded
as Banach algebras of functions on $E\cross E$ by Proposition 
\ref{prop:ehtpisalg}.  However, as pointed out in \cite{shulmant},
an element $u$ of
$\valb{E}$ may not be continuous on $E\cross E$, even if $E$ is compact.  
Hovever, if $\fC$ is a closed subalgebra of $\contb{E}$ (say
$\fC=\conto{E}$), then for each $u\in\fC\ehtens\fC\subset\valb{E}$, 
the pointwise slices, $u(\cdot,x)$ and $u(x,\cdot)$ for any fixed $x\iin E$,
will always be elements of $\fC$. In the case where
$E$ is a compact group, $\valoc{E}$ is discussed in \cite{spronkt},
and in a profound way in \cite{varopoulos}.  We note that by
Grothendieck's Inequality, $\valoc{E}=\conto{E}\tens^\gam\conto{E}$
(projective tensor product), up to equivalent norms.

If $u:E\cross E\to\Cee$, we say that $u$ is {\it positive definite} if
for any finite collection of elements $x_1,\dots,x_n$ from $E$, the matrix
$\sbrac{u(x_i,x_j)}$ is of positive type.

If $\fA$ is any abelian C*-algebra for which there is
a locally compact space $E$ and an injective $*$-homomorphism
$F:\fA\to\contb{E}$, then there is an isometric algebra homomorphism
$F\tens F:\fA\ehtens\fA\to\valb{E}$, 
by \cite{effrosr} or \cite[Cor.\ 2.3]{spronk}. 

The following theorem generalises \cite[Thm.\ 5.1]{stormer}.

\begin{posabelian}\label{theo:posabelian}
Let $\fA$ be an abelian C*-subalgebra of $\bdophil$ for which there
is a locally compact space $E$ and an injective $*$-homomorphism 
$F:\fA\to\contb{E}$.
If $T\in\fA\ehtens\fA$ and  $u=(F\tens F)T$, so $u\in F(\fA)\ehtens F(\fA) 
\subset\valb{E}$, then the following are equivalent:

{\bf (i)} $T$ is positive.

{\bf (ii)} $T$ is completely positive.

{\bf (iii)} $u$ is positive definite.
\end{posabelian}

\proof {\bf (i)}$\implies${\bf (ii)}  If $\fB$ is any maximal abelian
subalgebra of $\bdophil$ which contains $\fA$, then $T$ is
a $\fB$-bimodule map.  The result then follows from Lemma \ref{lem:posiscp}.

{\bf (ii)}$\implies${\bf (iii)} By Lemma \ref{lem:cpmapchar} we have that
$T=\sum_{i\in I}a_i\tens a_i^*$ for some family of elements from $\fA$ for 
which $\sum_{i\in I}a_ia_i^*$ converges weak* in $\bdophil$.
Let $\vphi_i=F(a_i)\iin\contb{E}$, so
\[
u=\sum_{i\in I}\vphi_i\tens\bar{\vphi}_i\quad\aand\quad
\unorm{\sum_{i\in I}|\vphi_i|^2}<\infty.
\]
Let $\xi:E\to\ltwo{I}$ be given by $\xi(x)=\bigl(\vphi_i(x)\bigr)_{i\in I}$.
Then for each $(x,y)\iin E\cross E$, we have that
\begin{equation}\label{eq:posdefchar}
u(x,y)=\inprod{\xi(x)}{\xi(y)}
\end{equation}
and hence $u$ is positive definite.

{\bf (iii)}$\implies${\bf (i)}  Since $u$ is positive definite function, then 
by \cite[\S 3.1]{guichardet}, there is a Hilbert space $\fL$ and a bounded 
function $\xi:E\to\fL$ such that (\ref{eq:posdefchar}) holds.  Let $p$ be
the orthogonal projection on $\fL$ whose range is 
$\wbar{\spn}\{\xi(x)\}_{x\in E}$, and let $\{\xi_i\}_{i\in I}$ be an 
orthonormal basis for $p\fL$.  Then for each $i$ the function
\[
\vphi_i=\inprod{\xi(\cdot)}{\xi_i}
\]
is in $F(\fA)$.  Indeed, given $\eps>0$ we can find $\alp_1,\dots,
\alp_n$ from $\Cee$ and $y_1,\dots,y_n$ from $E$, such that
\[
\norm{\xi_i-\sum_{k=1}^n\alp_k\xi(y_k)}<\eps
\]
whence
\[
\unorm{\vphi_i-\sum_{k=1}^n\bar{\alp}_k u(\cdot,y_k)}
=\unorm{\inprod{\xi(\cdot)}{\xi_i}-
\sum_{k=1}^n\bar{\alp}_k\inprod{\xi(\cdot)}{\xi(y_k)}}
<\unorm{\xi}\eps.
\]
Hence $\vphi_i$ can be uniformly approximated arbitrarily closely by elements 
of $F(\fA)$, and our conclusion holds.  It then follows by Parseval's 
Identity that for any $(x,y)\iin E\cross E$
\[
u(x,y)=\inprod{p\xi(x)}{p\xi(y)}
=\sum_{i\in I}\inprod{\xi(x)}{\xi_i}\inprod{\xi_i}{\xi(y)}
=\sum_{i\in I}\vphi_i(x)\wbar{\vphi_i(y)}.
\]
Hence we may write
\[
u=\sum_{i\in I}\vphi_i\tens\bar{\vphi}_i\quad\wwith\quad
\unorm{\sum_{i\in I}|\vphi_i|^2}=\unorm{\xi}^2<\infty.
\]
Letting $a_i=F^{-1}(\vphi_i)\iin\fA$, we get that 
$T=(F\tens F)^{-1}u=\sum_{i\in I}a_i\tens a_i^*$ and is thus positive. \endpf

\section{Representations of Groups in Completely Bounded Maps}

Let $G$ be a locally compact group, let $\fA$ be a unital Banach algebra 
which is also a dual space with predual $\fA_*$, 
and let $\alp:G\to\fA_\inv$ be a weak* continuous
bounded homomorphism where $\fA_\inv$ denotes the group of invertible
elements in $\fA$.  In particular we assume $\alp(e)$ is the unit of $\fA$
Denote the space of bounded complex Borel measures on 
$G$ by $\measg$.
Recall that $\measg$ is the dual space to the space $\conto{G}$ of
continuous functions vanishing at infinity.  Recall too that
$\measg$ is a Banach algebra via {\it convolution}: for each 
$\mu,\nu\iin\measg$ we define $\mu\con\nu$ by 
\begin{equation}\label{eq:convolution}
\int_G \vphi d\mu\con\nu=\int_G\int_G\vphi(st)d\mu(s)d\nu(t)
\end{equation}
for each $\vphi\iin\conto{G}$.  We note that since each of $\mu$ and
$\nu$ can be approximated in norm by compactly supported bounded measures,
(\ref{eq:convolution}) holds for any choice of $\vphi\iin\contb{G}$ too.
If $\mu\in\measg$, let
\[
\alp_1(\mu)=\text{weak*-}\int_G \alp(s)d\mu(s)
\]
i.e.\ if $\ome\in\fA_*$, then $\dpair{\alp_1(\mu)}{\ome}=
\int_G \dpair{\alp(s)}{\ome}d\mu(s)$.  Then $\alp_1:\measg\to\fA$
is a bounded linear map for if $\unorm{\alp}=\sup_{s\in G}\norm{\alp(s)}$, then
\begin{equation}\label{eq:repcontract}
\norm{\alp_1(\mu)}
=\sup_{\ome\in\ball{\fA_*}}\abs{\int_G \dpair{\alp(s)}{\ome}d\mu(s)}
\leq \int_G \unorm{\alp} d|\mu|(s)=\unorm{\alp}\norm{\mu}_1.
\end{equation}

Recall that the dual $\fA^*$ is a contractive $\fA$-bimodule
where for $b\iin\fA$ and $F\iin\fA^*$ we define $b\mult F$ and
$F\mult b$ in $\fA^*$ by $\dpair{a}{b\mult F}=\dpair{ab}{F}$ and
$\dpair{a}{F\mult b}=\dpair{ba}{F}$, for each $a\iin\fA$.  We say that
a subspace $\Omega$ of $\fA^*$ is a right $\alp(G)$-submodule if
$\ome\mult\alp(s)\in\Omega$, for each $\ome\iin\Omega$
and $s\iin G$.

\begin{alpishomo}\label{prop:alpishomo}
Let $G$, $\fA$ and $\alp$ be as above.  Moreover, suppose that
$\fA_*$ is both a left $\fA$-submodule of $\fA^*$ and 
a right $\alp(G)$-submodule.  Then $\alp_1:\measg\to\fA$
is a unital algebra homomorphism.
\end{alpishomo}

\proof If $\mu,\nu\in\measg$ and $\ome\in\fA^*$ then
\begin{align*}
\dpair{\alp_1(\mu)\alp_1(\nu)}{\ome}
&=\dpair{\alp_1(\mu)}{\alp_1(\nu)\mult\ome} \\
&=\int_G \dpair{\alp(s)}{\alp_1(\nu)\mult\ome}d\mu(s) \\
&=\int_G \dpair{\alp_1(\nu)}{\ome\mult\alp(s)}d\mu(s) \\
&=\int_G \int_G \dpair{\alp(t)}{\ome\mult\alp(s)}d\nu(t)d\mu(s) \\
&=\int_G \int_G \dpair{\alp(st)}{\ome}d\nu(t)d\mu(s).
\end{align*}
where the hypotheses guarantee that $\alp_1(\nu)\mult\ome\in\fA_*$
and that $\ome\mult\alp(s)\in\fA_*$, for each $s$.  
By Fubini's Theorem we have that
\[
\int_G \int_G \dpair{\alp(st)}{\ome}d\nu(t)d\mu(s)
=\int_G \int_G \dpair{\alp(st)}{\ome}d\mu(s)d\nu(t)
=\dpair{\alp_1(\mu\con\nu)}{\ome}
\]
where we 
note that $(s,t)\mapsto\dpair{\alp(st)}{\ome}$ is continuous and bounded,
hence $\mu\cross\nu$-integrable.  

That $\alp_1(\del_e)=\alp(e)$ follows from that $\fA_*$ is a separating
for $\fA$.  Hence $\alp_1$ is a unital map.
\endpf 

By a symmetric argument, the above proposition also holds if $\fA_*$ is
assumed to be both a right $\fA$-submodule of $\fA^*$ and a left
$\alp(G)$-submodule.

\begin{moduleexample}\label{ex:moduleexample}{\rm
{\bf (i)} Let $\fX$ be a Banach space admitting a predual $\fX_*$.
Then we have that $\fA=\bdop{\fX}$ is a dual unital Banach algebra admitting 
a predual $\fA_*=\fX\tens^\gam\fX_*$, via the dual pairing
\[
\dpair{T}{x\tens\ome}=\dpair{Tx}{\ome}\ffor T\iin\fA,\; x\iin\fX\aand
\ome\iin\fX_*.
\]
Here $\tens^\gam$ denotes the {\it projective tensor product}.  We have
then that $\fA_*$ is a left $\fA$ submodule of $\fA^*$.  Indeed, 
for any $S,T\iin\fA$ and elementary tensor $x\tens\ome\iin\fA_*$ we have that,
\[
\dpair{ST}{x\tens\ome}=\dpair{STx}{\ome}=\dpair{S}{(Tx)\tens\ome}
\]
so $T\mult (x\tens\ome)=(Tx)\tens\ome$.  

If $\wbdop{\fX}$ denotes the 
weak*-weak* continuous bounded linear maps on $\fX$ then 
$\fA_*$ is a right $\wbdop{\fX}$-submodule of $\fA^*$.  Thus we obtain the 
situation of Proposition \ref{prop:alpishomo} whenever 
$\alp:G\to\fA_\inv$ is a weak* continuous bounded homomorphism
whose range is in $\wbdop{\fX}$.  In particular, this happens
when $\fX$ is reflexive and $\alp$ is a non-degenerate
strong operator continuous representation on $\fX$.

{\bf (ii)}  The example above can be easily modified for the case
where $\fV$ is a dual operator space and $\fA=\cbop{\fV}$.

{\bf (iii)} There is a standard identification $\wcbop{\bdop{\hil}}
\cong\cbop{\fK(\hil),\bdop{\hil}}$, and thus an identification
of $\bdop{\hil}\ehtens\bdop{\hil}\cong\cbop{\fK(\hil),\bdop{\hil}}$.
In fact, as shown in \cite{blechers}, this latter identification is
a weak* homeomorphism.  Indeed, using standard identifications
with row and column Hilbert spaces and the {\it operator projective 
tensor product}, $\hat{\tens}$ (see \cite{blecherp} or
\cite[II.9.3]{effrosrB}), we have
\begin{align*}
\bdop{\hil}_*\htens\bdop{\hil}_*&\cong
\brac{\wbar{\hil}_r\htens\hil_c}\htens\brac{\wbar{\hil}_r\htens\hil_c} \\
&\cong\wbar{\hil}_r\htens\brac{\hil_c\htens\wbar{\hil}_r}\htens\hil_c 
\cong\wbar{\hil}_r\hat{\tens}\fK(\hil)\hat{\tens}\hil_c \\
&\cong\fK(\hil)\hat{\tens}\wbar{\hil}_r\hat{\tens}\hil_c
\cong\fK(\hil)\hat{\tens}\bdop{\hil}_*.
\end{align*}
On elementary tensors this identification is given by
\[
(\xi^*\tens\eta)\tens(\zeta^*\tens\vth)\mapsto
(\eta\tens\zeta^*)\tens(\xi^*\tens\vth)
\]
where for vectors $\xi,\eta\iin\hil$ we let $\xi\tens\eta^*$ denote the
usual rank 1 operator and $\xi^*\tens\eta$ the usual vector functional.
Now if $T=\sum_{i\in I}a_i\tens b_i\iin\bdop{\hil}\ehtens\bdop{\hil}$
then, in the dual pairing (\ref{eq:htpdual}), we have that
\[
\dpair{T}{(\xi^*\tens\eta)\tens(\zeta^*\tens\vth)}
=\sum_{i\in I}\inprod{a_i\eta}{\xi}\inprod{b_i\vth}{\zeta}.
\]
Meanwhile, in the 
$\cbop{\fK(\hil),\bdop{\hil}}$--$\fK(\hil)\hat{\tens}\bdop{\hil}_*$ 
duality we have that
\begin{align*}
\dpair{T}{(\eta\tens\zeta^*)\tens(\xi^*\tens\vth)}
&=\dpair{\sum_{i\in I} a_i\eta\tens(b_i^*\zeta)^*}{\xi^*\tens\vth} \\
&=\sum_{i\in I}\inprod{a_i\eta}{\xi}\inprod{b_i\vth}{\zeta}.
\end{align*}

Now for every elementary tensor $k\tens\ome
\iin\fK(\hil)\hat{\tens}\bdop{\hil}_*$ and 
$T\iin\bdop{\hil}\ehtens\bdop{\hil}$, we have that 
$(k\tens\ome)\mult T=k\tens(\ome\mult T)\in\fK(\hil)\hat{\tens}\bdop{\hil}_*$.
Hence $\bdop{\hil}_*\htens\bdop{\hil}_*\cong\fK(\hil)\hat{\tens}\bdop{\hil}_*$
is a right module for $\bdop{\hil}\ehtens\bdop{\hil}$.  We note that
$\bdop{\hil}_*\htens\bdop{\hil}_*$ is a left 
$\bdop{\hil}\htens\bdop{\hil}$-module. \endpf
}\end{moduleexample}

We will identify the group algebra $\bloneg$ with the closed ideal in
$\measg$ of measures which are absolutely continuous with respect to the 
left Haar measure $m$ (whose integral we will denote $\int_G\cdots ds$).
We will identify the discrete group algebra $\loneg$ with the closed
subspace of $\measg$ generated by all of the Dirac measures 
$\{\del_s:s\in G\}$.  We let
\[
\measalp=\wbar{\alp_1(\measg)},\quad
\cmeasalp=\wbar{\alp_1(\bloneg)}\quad\aand\quad
\dmeasalp=\wbar{\alp_1(\loneg)}
\]
where each of the closures is in the norm topology of $\fA$.

The following proposition is surely well-known, though we have been unable
to find it in the literature.

\begin{repwunit}\label{prop:repwunit}
Given $G$, $\fA$ and $\alp$ satifying the hypotheses of Proposition
\ref{prop:alpishomo}, $\alp$ is norm continuous
if and only if $\cmeasalp$ is unital.
\end{repwunit}

\proof Let $(e_U)$ be the bounded approximate
identity for $\bloneg$ given by $e_U=\frac{1}{m(U)}1_U$ (normalised indicator
function), indexed over the family of all relatively compact neighbourhoods of 
the identity $e$ in $G$, partially ordered by reverse inclusion.

$\rif$ Let $\eps>0$. Let $V$ be any relatively compact neighbourhood 
of $e$ for which $\norm{\alp(s)-\alp(e)}<\eps$ for each $s\iin V$. 
Then for any relatively compact neighbourhood
$U$ of $e$ which is contained in $V$ we have
\begin{align*}
\norm{\alp_1(e_U)-\alp(e)}&=\norm{\frac{1}{m(U)}\int_U \alp(s)ds-\alp(e)} \\
&\leq\frac{1}{m(U)}\int_U\norm{\alp(s)-\alp(e)}ds
<\eps
\end{align*}
where the second from last inequality is proved just as in 
(\ref{eq:repcontract}).
Thus $\alp(e)=\lim_U\alp_1(e_U)$ in norm, so $\alp(e)\in\cmeasalp$. 
Now $\alp(e)$ is the unit for $\fA$, and hence the unit for $\cmeasalp$.

$\lif$ It is a standard fact that $\lim_U\alp_1(e_U)=\alp(e)$ in the weak* 
topology of $\fA$.  Indeed, $\lim_U \int_G e_U(s)\vphi(s)ds=\vphi(e)$ for any
continuous function $\vphi$; set $\vphi=\dpair{\alp(\cdot)}{\ome}$ for
any $\ome\iin\fA_*$.  Now let $E$ be the unit for $\cmeasalp$.  We will
establish that $E=\alp(e)$, the unit of $\fA$.  
First, the map $s\mapsto\alp(s)E$ is norm 
continuous.  Indeed $E\in\cmeasalp$ and can thus be norm approximated by
$\{\alp_1(f):f\in\bloneg\}$.  Moreover, if $f\in\bloneg$ then we have that
\[
\norm{\alp(s)\alp_1(f)-\alp_1(f)}
=\norm{\alp_1(\del_s\con f-f)} 
\leq\unorm{\alp}\norm{\del_s\con f-f}_1
\overset{s\to e}{\longrightarrow}0
\]
where the inequality follows from (\ref{eq:repcontract}) and 
limit follows from \cite[20.4]{hewittrI}.
Next, for any compact neighbourhood $U$ of $e$ we have that
\[
\alp_1(e_U)=\alp_1(e_U)E=\frac{1}{m(U)}\int_U\alp(s)ds\mult E
=\frac{1}{m(U)}\int_U\alp(s)Eds
\]
where we note that right multiplication is weak*-continuous in $\fA$, by 
hypothesis.
Now, let $\eps>0$ be given, and find a neighbourhood $V$ of $e$ in $G$
such that $\norm{\alp(s)E-E}<\eps$ for each $s\iin V$.  Then for any
relatively compact neighbourhood $U$ of $E$, contained in $V$, we have that
\[
\norm{\alp_1(e_U)-E}=\norm{\frac{1}{m(U)}\int_U \alp(s)Eds-E}
\leq\frac{1}{m(U)}\int_U\norm{\alp(s)E-E}ds<\eps
\]
where the second from last inequality is proved just as in 
(\ref{eq:repcontract}).  Hence we have that $\lim_U\alp_1(e_U)=E$ in norm, so,
{\it a fortiori}, weak$^*$-$\lim_U\alp_1(e_U)=E$.  It then follows from
above that $E=\alp(e)$, so $\alp(e)\in\cmeasalp$.  Thus 
\[
\norm{\alp(s)-\alp(e)}=\norm{\alp(s)E-E}\overset{s\to e}{\longrightarrow}0.
\]
Hence $\alp$ is norm continuous at $e$, and thus norm continuous on all of $G$.
\endpf

\begin{repwunit1}\label{cor:repwunit1}
For $G$, $\fA$ and $\alp$ as above, the following are equivalent:

{\bf (i)} $\alp$ is norm continuous \quad\quad
{\bf (ii)} $\cmeasalp=\measalp$ \quad\quad
{\bf (iii)} $\cmeasalp=\dmeasalp$.
\end{repwunit1}

\proof (i)$\iff$(ii) If $\alp$ is norm continuous, then $\cmeasalp$ contains
the unit $\alp(e)$ by Proposition \ref{prop:repwunit}.  Hence,
$\cmeasalp$ is an ideal in $\measalp$, containing the unit.
Conversely, if $\cmeasalp=\measalp$ then $\cmeasalp$ is unital,
and norm continuity of $\alp$ follows from Proposition \ref{prop:repwunit}.

(i)$\implies$(iii) Since (ii) holds, the inclusion
$\cmeasalp\supset\dmeasalp$ is clear.  To obtain the opposite inclusion, 
note that for any continuous function of compact support $\vphi$ 
-- the family of which 
is dense in $\bloneg$ -- the function $s\mapsto\vphi(s)\alp(s)$, 
from $G$ to $\dmeasalp$, can be uniformly approximated by Borel simple
functions.  
Hence $\alp_1(\vphi)=\int_G\vphi(s)\alp(s)ds$
may be regarded as a Bochner integral, and is thus in $\dmeasalp$,
since each $\alp(s)\in\dmeasalp$.
It then follows that $\alp_1(\bloneg)\subset\dmeasalp$ and hence 
$\cmeasalp\subset\dmeasalp$.

(iii)$\implies$(i)  Since $\cmeasalp\supset\dmeasalp$, $\cmeasalp$ is
unital, and the result follows from Proposition \ref{prop:repwunit}. \endpf

Now suppose that $\pi:G\to\unophil\subset\bdophil_\inv$ is a strongly 
continuous unitary representation (which is equivalent to it being 
weak* continuous).  We will define $\pi_1:\measg\to\bdophil$ as above,
but will use the notation 
\[
\mstarpi=\wbar{\pi_1(\measg)},\quad 
\cstarpi=\wbar{\pi_1(\bloneg)}\quad\aand\quad
\cstardpi=\wbar{\pi_1(\loneg)}
\]
to indicate that these are C*-algebras.  Using von Neumann's
double commutant theorem, we have that $\cstarpi$
and $\cstardpi$ each generate the same von Neumann algebra,
$\vnpi$.  We note that $\mstarpi\subset\vnpi$ but $\mstarpi\not=\vnpi$
in general.  Thus, in particular, there is no reason to suspect that 
$\mstarpi$ is a dual space.

\begin{repsingeneral}\label{prop:repsingeneral}
If $\mu\in\measg$, then
\begin{equation}\label{eq:gammamap}
\gampi(\mu)=\int_G \pi(s)\tens\pi(s)^*d\mu(s)
\end{equation}
defines an element of  $\bdophil\whtens\bdophil$, and the integral converges
in the weak* topology, i.e.\ for each $x\iin\bdophil_*\htens\bdophil_*$,
\[
\dpair{\gampi(\mu)}{x}=\int_G\dpair{\pi(s)\tens\pi(s)^*}{x}d\mu(s).
\]
Moreover, 

{\bf (i)} $\gampi:\measg\to\bdophil\whtens\bdophil$ is a contractive
homomorphism whose range is contained in the algebra $\mstarpi\ehtens\mstarpi$.

{\bf (ii)} $\gampi(\bloneg)\subset\cstarpi\ehtens\cstarpi$.

{\bf (iii)} $\gampi(\loneg)\subset\cstardpi\htens\cstardpi$.

{\bf (iv)} If $\pi$ is norm continuous, then $\gampi(\measg)
\subset\cstardpi\htens\cstardpi$.
\end{repsingeneral}

\proof 
{\bf (i)} First, let us see that, for each $\mu\iin\measg$, the integral
in (\ref{eq:gammamap}) converges as claimed.  
This amounts to verifying that $s\mapsto\pi(s)\tens\pi(s)^*$ is a 
weak* continuous representation from $G$ into $(\bdophil\whtens\bdophil)_\inv$,
i.e.\ that $s\mapsto\dpair{\pi(s)\tens\pi(s)^*}{x}$ 
is continuous for each
$x\iin\bdophil_*\htens\bdophil_*$, by (\ref{eq:htpdual}).  
If $x\in\bdophil_*\htens\bdophil_*$
and $\eps>0$, then there is $x_\eps\iin\bdophil_*\tens\bdophil_*$ such that
$\hnorm{x-x_\eps}<\eps$.  The function $x_{\eps,\pi}$, given by
$s\mapsto\dpair{\pi(s)\tens\pi(s)^*}{x_\eps}$, is clearly 
continuous on $G$, and $\unorm{x_\pi-x_{\eps,\pi}}\leq\hnorm{x-x_\eps}
<\eps$.  Thus, taking choices of $\eps$ tending to $0$, we see that
$x_\pi$ is a continuous function on $G$.

Since $\whnorm{\pi(s)\tens\pi(s)^*}=1$ for each $s\iin G$, the contractivity 
of $\gampi$ follows from 
(\ref{eq:repcontract}).  That $\gampi$ is a homomorphism 
follows from Proposition \ref{prop:alpishomo} and Example 
\ref{ex:moduleexample} (iii). 

To see that $\gampi(\mu)\in\mstarpi\ehtens\mstarpi$, for any given $\mu\iin
\measg$, we will inspect the image of a typical weak*-weak* continuous left 
slice map on $\gampi(\mu)$ and use \cite[Thm.\ 2.2]{spronk}. 
If $\ome\in\bdophil_*$, then
\begin{equation}\label{eq:slicetomstar}
L_\ome\left(\gampi(\mu)\right)
=\int_G \dpair{\pi(s)}{\ome}\pi(s)^*d\mu(s)
=\int_G \pi(s)d(\ome_\pi\mu)^\vee(s)\in\mstarpi
\end{equation}
where $\ome_\pi\mu$ is the measure with Radon derivative
$d(\ome_\pi\mu)/d\mu=\ome_\pi$ (here $\ome_\pi(s)=\dpair{\pi(s)}{\ome}$), 
and $\nu^\vee(E)=\nu(E^{-1})=\wbar{\nu^*(E)}$
for any Borel measure $\nu$.  The computation for any right slice map is 
similar.

{\bf (ii)} This follows from a computation similar to (\ref{eq:slicetomstar}).

{\bf (iii)} If $\mu=\sum_{s\in G}\alp(s)\del_s$, where 
$\sum_{s\in G}|\alp(s)|<\infty$, then since $\pi(s)\tens\pi(s)^*
\in\cstardpi\htens\cstardpi$ for each $s\iin G$, it follows too that
\[
\gampi(\mu)=\sum_{s\in G}\alp(s)\pi(s)\tens\pi(s)^*\in\cstardpi\htens\cstardpi.
\]

{\bf (iv)}  If we let $\alp:G\to(\bdophil\whtens\bdophil)_\inv$ be
given by $\alp(s)=\pi(s)\tens\pi(s)^*$, then $\alp$ is norm continuous.
Hence
\[
\gampi(\measg)\subset\measalp=\dmeasalp\subset\cstardpi\htens\cstardpi
\]
by (iii) above and Corollary \ref{cor:repwunit1}.  \endpf

\begin{contincb}\label{rem:contincb}{\rm
We note that (\ref{eq:gammamap}) also converges in the 
$\cbop{\bdophil}$--$(\bdophil\what{\tens}$ $\bdophil_*)$ topology.  Indeed, if
$a\in\bdophil$ and $\eta^*\tens\xi$ is a vector functional in $\bdophil_*$,
then for any $s\iin G$ we have that
\[
\dpair{\pi(s)\tens\pi(s)^*}{a\tens(\eta^*\tens\xi)}
=\dpair{\pi(s)a\pi(s)^*}{\eta^*\tens\xi}
=\dpair{a}{(\pi(s)^*\eta)^*\tens\pi(s)^*\xi}
\]
where $s\mapsto(\pi(s)^*\eta)^*\tens\pi(s)^*\xi$ is continuous in the
norm topology of $\bdophil_*$.  Hence $s\mapsto
\dpair{\pi(s)\tens\pi(s)^*}{a\tens(\eta^*\tens\xi)}$ is continuous.
In particular, for each $a\iin\bdophil$ and $\mu\iin\measg$ we have that
\[
\gampi(\mu)a=\int_G\pi(s)a\pi(s)^*d\mu(s)
\]
where the integral converges in the weak* topology of $\bdophil$

We observe that it is possible, for each $\mu\iin\measg$, to see  
$\gampi(\mu)|_{\fK(\hil)}$ as an integral converging in the point-norm 
topology.  However, our approach for obtaining (\ref{eq:gammamap}) better
lends itself to (\ref{eq:pointeval}).} \endpf
\end{contincb}

We let the {\it augmentation ideal} in $\bloneg$ be given by
\[
\augidealg=\cbrac{f\in\bloneg:\int_Gf(s)ds=0}.
\]

\begin{repinhtp}\label{theo:repinhtp}
For any strongly continuous representation $\pi:G\to\unophil$, the following
are equivalent:

{\bf (i)} $\pi$ is norm continuous.

{\bf (ii)} $\gampi(\bloneg)\subset\cstarpi\htens\cstarpi$.

{\bf (iii)} there is an $f\iin\bloneg\setdif\augidealg$ such that
$\gampi(f)\in\cstarpi\htens\cstarpi$.
\end{repinhtp}

\proof That (i) implies (ii) follows from Proposition \ref{prop:repsingeneral}
(iv) and the fact that $\cstarpi=\cstardpi$.  That (ii) implies (iii) is
trivial.  Suppose now that $f$ satisfies statement (iii).  Without loss of 
generality, we may suppose that $\int_G f(s)ds=1$.  Then by
\cite{blechers}, there exist sequences $\{a_i\}_{i\in\En}$ and 
$\{b_i\}_{i\in\En}$ from $\cstarpi$ such that $\sum_{i=1}^\infty a_ia_i^*$
and $\sum_{i=1}^\infty b_i^*b_i$ converge in norm, and
\[
\gampi(f)x=\sum_{i=1}^\infty a_i x b_i
\]
for each $x\in \bdophil$.  But it then follows from Remark \ref{rem:contincb}
that
\[
I=\int_G f(s) \pi(s)I\pi(s)^*ds=\gampi(f)I=
\sum_{i=1}^\infty a_ib_i\in\cstarpi.
\]
Hence $\pi$ is norm continuous by Proposition \ref{prop:repwunit}.  \endpf

In the next section, we will address the necessity of the assumption that 
$f\in\bloneg\setdif\augidealg$ in (iii) above.

It is interesting to note that the kernel of $\gampi$ is related to the
kernel of a more familiar representation.  Below, we will let
$\wbar{\hil}$ denote the conjugate Hilbert space and 
$\bar{\pi}:G\to\unop{\wbar{\hil}}$ denote the conjugate representation.
We will also let $\pi\tens\bar{\pi}:G\to\unop{\hil\tens_2\wbar{\hil}}$
be the usual tensor product of representations on the Hilbert space
$\hil\tens_2\wbar{\hil}$.

\begin{kergampi}\label{prop:kergampi}
$\ker\gampi=\ker(\pi\tens\bar{\pi})_1$.
\end{kergampi}

\proof We have that $\mu\in\ker\gampi$ if and only if 
\[
0=\dpair{\gampi(\mu)}{\ome_{\xi,\eta}\tens\ome_{\zeta,\vth}}
\]
for every elementary tensor 
of vector functionals $\ome_{\xi,\eta}\tens\ome_{\zeta,\vth}$
in $\bdop{\hil}_*\htens\bdop{\hil}_*$.  (Note that we earlier
had used the notation $\ome_{\xi,\eta}=\eta^*\tens\xi$.)
We may compute
\begin{align*}
\dpair{\gampi(\mu)}{\ome_{\xi,\eta}\tens\ome_{\zeta,\vth}}
&=\int_G\inprod{\pi(s)\xi}{\eta}\inprod{\pi(s)^*\zeta}{\vth}d\mu(s) \\
&=\int_G\inprod{\pi(s)\xi}{\eta}\wbar{\inprod{\pi(s)\vth}{\zeta}}d\mu(s) \\
&=\int_G\inprod{\pi\tens\bar{\pi}(s)\;\xi\tens\bar{\vth}}
               {\eta\tens\bar{\zeta}}d\mu(s) \\
&=\inprod{(\pi\tens\bar{\pi})_1(\mu)\;\xi\tens\bar{\vth}}
               {\eta\tens\bar{\zeta}}.
\end{align*}
Thus it follows that $\mu\in\ker\gampi$ if and only if 
$\mu\in\ker(\pi\tens\bar{\pi})_1$. \endpf

In particular, if we let $\coef{\pi\tens\bar{\pi}}$ be the linear space
generated by all of the coefficient functions, $s\mapsto
\inprod{\pi\tens\bar{\pi}(s)\;\xi\tens\bar{\vth}}{\eta\tens\bar{\zeta}}$,
we see that $\mu\in\ker\gampi$ exactly when $\mu$,  as a functional on 
$\contb{G}$, annihilates $\coef{\pi\tens\bar{\pi}}$.

\section{Abelian Groups}
\label{sec:abelian}

For this section we let $G$ be a locally compact {\it abelian} group, and
we let $\dgrp$ denote its topological dual group.  For each $s\iin G$, we 
will let $\hat{s}$ denote the associated unitary character on $\dgrp$, defined
by $\hat{s}(\sig)=\sig(s)$ for each $\sig\iin\dgrp$. 

As above, we will let
$\pi:G\to\unophil$ be a strongly continuous unitary representation.
We let $E_\pi$ denote the spectrum of $\cstarpi$.  Since $\cstarpi$ is
a quotient of the enveloping C*-algebra $\cstarg$, and 
$\cstarg\cong\conto{\dgrp}$, we may consider $E_\pi$ to be a closed subset 
of $\dgrp$. Moreover, the natural isomorphism $F_\pi:\cstarpi\to\conto{E_\pi}$
satisfies 
\[
F_\pi\bigl(\pi_1(f)\bigr)=\hat{f}|_{E_\pi}
\]
for each $f\iin\bloneg$, where $\hat{f}(\sig)=\sig_1(f)=\int_G f(s)\sig(s)ds$
for each $\sig\iin\dgrp$.  We note that our notation $f\mapsto\hat{f}$,
for the Fourier transform, differs from that of our main reference,
\cite{hewittrI}.  

We would like to be able to extend $F_\pi$ to some suitable map
$\bar{F}_\pi$ on $\vnpi$.  It is not clear that this can be done in general,
but it can be done in many cases.

\begin{abelianreps0}\label{lem:abelianreps0}
Consider the following conditions for $\pi$ or $G$ below:

{\bf (a)} $\hil$ admits a maximal countable family of mutually orthogonal 
cyclic subspaces for $\pi$.

{\bf (b)} There is a family $\{U_i\}_{i\in I}$ of separable open
subsets of $E_\pi$ such that $E_\pi=\dot{\bigcup}_{i\in I}U_i$.

{\bf (c)} $\dgrp$ has a separable open subgroup.

{\bf (d)} $G$ is compactly generated.

\noindent Then, under any one of these conditions there exists a regular 
Borel measure $\nu$ on $E_\pi$,
bounded on compacta, such that there is a normal $*$-homomorphism
$\bar{F}_\pi:\vnpi\to\blinfty{E_\pi,\nu}$ which extends $F_\pi$.
\end{abelianreps0}

\proof {\bf (a)} By standard arguments (see \cite[\S 7]{dixmier}, 
for example), $\vnpi$ admits a faithful normal
state $\ome$.  Then the measure $\nu$ given by
\begin{equation}\label{eq:measdef}
\int_{E_\pi}\vphi(\sig)d\nu(\sig)=\ome(F_\pi^{-1}\vphi)
\end{equation}
for each $\vphi\iin\conto{E_\pi}$, gives rise to the desired map $\bar{F}_\pi$.

{\bf (b)} Since $E_\pi=\dot{\bigcup}_{i\in I}U_i$, we have that
$\conto{E_\pi}=\ceeos\bigoplus_{i\in I}\conto{U_i}$.  If we let 
$\fC_i=F_\pi^{-1}\bigl(\conto{U_i}\bigr)$, then 
$\fM_i=\wbar{\fC_i}^{w^*}$ is an ideal in $\vnpi$.  
The ideals $\fM_i$ are mutually orthogonal, and hence if $\{p_i\}_{i\in I}$
is the family of projections for which $\fM_i=p_i\vnpi$ for each $i$,
then $\sum_{i\in I}p_i=I$. Since each $\fC_i$ is separable, each $\fM_i$
is countably generated, and hence there is a normal state $\ome_i$ on $\vnpi$
with support projection $p_i$.  Let $\nu_i$ be the measure on $E_\pi$
associated with $\ome_i$ as in (\ref{eq:measdef}).  Then $\supp{\nu_i}=U_i$
for each $i$, and $\nu=\bigoplus_{i\in I}\nu_i$ is the desired measure.

{\bf (c)} If $\dgrp$ has a separable open subgroup $X$, let
$T$ be any transversal for $X$ in $\dgrp$, and we have that
$E_\pi=\dot{\bigcup}_{\tau\in T}\brac{E_\pi\cap\tau X}$, and again
we obtain (b).

{\bf (d)} If $G$ is compactly generated, then by \cite[9.8]{hewittrI}
there is a topological isomorphism $G\cong\Zee^n\cross\Ree^m\cross K$,
where $K$ is compact.  Then $\dgrp\cong\Tee^n\cross\Ree^m\cross\what{K}$,
and the subgroup $X$ corresponding to $\Tee^n\cross\Ree^m$ is open and 
separable, and hence (c) holds.  \endpf

We will need to use an extension of $F_\pi$ of a different
nature than in the lemma above. Since $\cstarpi$ is an essential ideal in 
$\mstarpi$, the map $F_\pi:\cstarpi\to\conto{E_\pi}$ extends to an injective
$*$-homomorphism $\til{F}_\pi:\mstarpi\to\contb{E_\pi}$, such that
$F_\pi(na)=\til{F}_\pi(n)F_\pi(a)$ for each $n\iin\mstarpi$
and $a\iin\cstarpi$, by
\cite[3.12.8]{pedersen}.  We note that  for each $\mu\iin\measg$,
\begin{equation}\label{eq:tilfpi}
\til{F}_\pi(\pi_1(\mu))=\hat{\mu}|_{E_\pi}
\end{equation}
where for each $\mu\iin\measg$, $\hat{\mu}(\sig)=\sig_1(\mu)
=\int_G\sig(s)d\mu(s)$.  Thus $\mu\mapsto\hat{\mu}$ is the Fourier-Stieltjes
transform.  To see the validity of (\ref{eq:tilfpi}), observe that
for each $f\iin\bloneg$ we have
\[
\hat{\mu}\hat{f}|_{E_\pi}=\what{\mu\con f}|_{E_\pi}
=F_\pi(\pi_1(\mu\con f))=\til{F}_\pi(\pi_1(\mu))F_\pi(f)=
\til{F}_\pi(\pi_1(\mu))\hat{f}|_{E_\pi}.
\]
Thus it follows that $\til{F}_\pi(\pi_1(\mu))\vphi=\hat{\mu}\vphi$ for each
$\vphi\iin\conto{E_\pi}$.

If any of the conditions of Lemma \ref{lem:abelianreps0} hold, 
then there exists a measure $\nu$ for which there is a normal extension 
$\bar{F}_\pi:\vnpi\to\blinfty{E_\pi,\nu}$ of $F_\pi$.  Then
for any $\mu\iin\measg$,
\begin{equation}\label{eq:baristilde}
\bar{F}_\pi(\pi_1(\mu))=\hat{\mu}|_{E_\pi}
\end{equation}
where we identify $\contb{E_\pi}$ as a closed subspace of 
$\blinfty{E_\pi,\nu}$.  To see (\ref{eq:baristilde}), we note that if
$(a_\beta)$ is any bounded approximate identity in $\cstarpi$, then
$\text{weak*-}\lim_\beta a_\beta=I$ in $\vnpi$, thus
$\text{weak*-}\lim_\beta F_\pi(a_\beta)=1_{E_\pi}$.  Hence
\begin{align*}
\bar{F}_\pi(\pi_1(\mu))
&=\text{weak*-}\lim_\beta \bar{F}_\pi(\pi_1(\mu)a_\beta)
=\text{weak*-}\lim_\beta F_\pi(\pi_1(\mu)a_\beta) \\
&=\text{weak*-}\lim_\beta \til{F}_\pi(\pi_1(\mu))F_\pi(a_\beta)
=\til{F}_\pi(\pi_1(\mu))=\hat{\mu}|_{E_\pi}.
\end{align*}

We will make use of the spaces $\valb{E}$, $\valo{E}$ and $\valoc{E}$, which
were defined in (\ref{eq:valgebras}).
If $\nu$ is any non-negative measure on $E$, we let
\[
\val{E,\nu}=\blinfty{E,\nu}\ehtens\blinfty{E,\nu}.
\]
Spaces of this type are discussed in \cite{spronk}.

\begin{abelianreps}\label{theo:abelianreps}
If $G$ is a locally compact abelian group and $\pi:G\to\unophil$
is a strongly continuous unitary representation, then for any $\mu$ in
$\measg$ and $(\sig,\tau)\iin E_\pi\cross E_\pi$ we have that
\[
(\til{F}_\pi\tens \til{F}_\pi)\gampi(\mu)(\sig,\tau)=\hat{\mu}(\sig\tau^{-1}).
\]
In particular, if $E$ is any closed subset of $\dgrp$ and
$\mu\in\measg$, then $u(\sig,\tau)=\hat{\mu}(\sig\tau^{-1})$ is
an element of $\valb{E}$, and $u\in\valo{E}$ if $\mu\in\bloneg$.  
Moreover, if $E$ is compact then $u\in\valoc{E}$.
\end{abelianreps}

\proof The result will be established in three stages.  The first two of 
these require additional hypotheses and are preparatory for the general case.

{\bf I.} Suppose that any one of the conditions of 
Lemma \ref{lem:abelianreps0} is satisfied.  Let $\nu$ be the measure
on $E_\pi$ and let $\bar{F}_\pi:\vnpi\to\blinfty{E_\pi,\nu}$
be the map given there.

If $\mu\in\measg$, we have that
\begin{align}\label{eq:pointeval}
(\til{F}_\pi\tens \til{F}_\pi)\gampi(\mu)&=
(\bar{F}_\pi\tens\bar{F}_\pi)\int_G \pi(s)\tens\pi(s)^*d\mu(s) \notag \\
&= \int_G \hat{s}|_{E_\pi}\tens\bar{\hat{s}}|_{E_\pi}d\mu(s) 
\end{align}
where the latter integral converges in the weak* topology of $\val{E_\pi,\nu}$.

For $(\sig,\tau)\iin E_\pi\cross E_\pi$ let
\[
u(\sig,\tau)=\hat{\mu}(\sig\tau^{-1}).
\]
Then $u\in\valb{E_\pi}$.  Indeed, we have that $\hat{\mu}\in\mathrm{B}(\dgrp)$,
the {\it Fourier-Stieltjes algebra} which is defined in \cite{eymard}.  
Thus there is a Hilbert space $\fL$, a continuous unitary representation 
$\rho:G\to\fU(\fL)$, and vectors $\xi,\eta\iin\fL$ with
$\norm{\mu}=\norm{\xi}\norm{\eta}$, such that $\hat{\mu}(\sig)=
\inprod{\rho(\sig)\xi}{\eta}$ for each $\sig\iin\dgrp$.  If 
$\{\xi_i\}_{i\in I}$ is an orthonormal basis for $\fL$, then we have,
using Parseval's formula, that 
\[
\hat{\mu}(\sig\tau^{-1})=\inprod{\rho(\sig\tau^{-1})\xi}{\eta}
=\sum_{i\in I}\inprod{\rho(\sig)\xi}{\xi_i}\inprod{\xi_i}{\rho(\tau)\eta}
\]
for any $(\sig,\tau)\iin E_\pi\cross E_\pi$.  Hence
\[
u=\sum_{i\in I}\inprod{\rho(\cdot)\xi}{\xi_i}\tens
\wbar{\inprod{\rho(\cdot)\eta}{\xi_i}}\in\valb{E_\pi}
\]
with $\ehnorm{u}\leq\norm{\xi}\norm{\eta}=\norm{\mu}$.
(This is similar to the proof of \cite[Prop.\ 5.1]{spronkt}.)
We note that if $\mu\in\bloneg$, then $\rho$ can be taken to
be the left regular representation and hence each 
$\inprod{\rho(\cdot)\xi}{\xi_i}$ and $\inprod{\rho(\cdot)\eta}{\xi_i}$
is in $\conto{E_\pi}$.  Hence, in this case we would have that 
$u\in\valo{E_\pi}$.

We wish to establish that
\begin{equation}\label{eq:uismu}
u=(\til{F}_\pi\tens \til{F}_\pi)\gampi(\mu).
\end{equation}
We will do this by using the dual pairing (\ref{eq:htpdual}).  
If $g\tens h$ is an elementary
tensor in $\blone{E_\pi,\nu}\htens\blone{E_\pi,\nu}$, then
\begin{align*}
\dpair{u}{g\tens h}
&=\int_{E_\pi}\int_{E_\pi} g(\sig)h(\tau)\hat{\mu}(\sig\tau^{-1})
  d\nu(\sig)d\nu(\tau) \\
&=\int_{E_\pi}\int_{E_\pi} g(\sig)h(\tau)
  \brac{\int_G \sig(s)\wbar{\tau(s)}d\mu(s)}d\nu(\sig)d\nu(\tau) \\
&=\int_G \brac{\int_{E_\pi}g(\sig)\hat{s}(\sig)d\nu(\sig)}
  \brac{\int_{E_\pi}h(\tau)\wbar{\hat{s}(\tau)}d\nu(\tau)}d\mu(s)
\end{align*}
where the version of Fubini's Theorem required is \cite[13.10]{hewittrI},
noting that $g$ and $h$ each have $\nu$-$\sig$-finite supports.
On the other hand, by (\ref{eq:pointeval}),
\begin{align*}
\langle(\til{F}_\pi\tens \til{F}_\pi)\gampi(\mu),&g\tens h\rangle
=\dpair{\int_G \hat{s}|_{E_\pi}\tens\bar{\hat{s}}|_{E_\pi}d\mu(s)}{g\tens h} \\
&=\int_G \brac{\int_{E_\pi}g(\sig)\hat{s}(\sig)d\nu(\sig)}
  \brac{\int_{E_\pi}h(\tau)\wbar{\hat{s}(\tau)}d\nu(\tau)}d\mu(s)
\end{align*}
and this shows that (\ref{eq:uismu}) holds.

{\bf II.} Suppose that $\mu$ is supported on a compactly generated open 
subgroup $H$ of $G$.

Let us first compute the spectrum $E_{\pi|_H}$ of $\cstarpih$.  
We note that $\dgrph=\dgrp|_H$ and that
the restriction map $r:\dgrp\to\dgrp|_H$ is a homomorphic topological
quotient map by \cite[24.5]{hewittrI}.  Moreover, $\ker r$ is
compact, by \cite[23.29(a)]{hewittrI}.
Then $E_{\pi|_H}=r(E_\pi)$.  To see this, 
observe that the map $\iota:\blone{H}\to\bloneg$, which we define to be
the inverse of $f\mapsto f|_H$,  extends
to an injective $*$-homomorphism $\iota_\pi:\cstarpih\to\cstarpi$.
In particular, then, each multiplicative linear functional on $\cstarpih$ is 
necessarily the restriction of such a functional on $\cstarpi$.
Let $r_\pi=r|_{E_\pi}$.
Then, the map $r_\pi:E_\pi\to r(E_\pi)$ induces an injective $*$-homomorphism
$j_{r_\pi}:\conto{r(E_\pi)}\to\conto{E_\pi}$, whose image is the subalgebra
of all functions which are constant on relative cosets of $\ker r$ in $E_\pi$. 
Now, if $g\in\blone{H}$, and $\sig\in E_\pi$ then
\begin{equation}\label{eq:rpiisquot}
\what{\iota g}(\sig)=\int_G \iota g(s)\sig(s)ds
=\int_H g(s)r(\sig)(s)ds=\hat{g}(r_\pi(\sig))
=j_{r_\pi}\hat{g}(\sig)
\end{equation}
from which it follows that every character on $\cstarpih$ is from
$r_\pi(E_\pi)$.  Moreover, it follows from (\ref{eq:rpiisquot}) that
\[
F_\pi\comp\iota_\pi=j_{r_\pi}\comp F_{\pi|_H}.
\]

Now, we let $\til{\iota}:\meash\to\measg$ be the homomorphism whose inverse
is $\kappa\mapsto\kappa_H$, where for any Borel subset $B$ of $G$, 
$\kappa_H(B)=\kappa(B\cap H)$.  Then $\iota$ induces an injective 
$*$-homomorphism $\til{\iota}_\pi:\mstarpih\to\mstarpi$.  It follows from the 
discussion above that $\tilfpih:\mstarpih\to\contb{r(E_\pi)}$.  Then
\begin{equation}\label{eq:jfpihisfpii}
\tiljrpi\comp \tilfpih =\til{F}_\pi\comp\til{\iota}_\pi
\end{equation}
where $\tiljrpi:\contb{E_{\pi|_H}}\to\contb{E_\pi}$ is the map
induced by $r_\pi:E_\pi\to E_{\pi|_H}$.  Indeed, if $\kappa\in\meash$,
then for each $\sig\iin E_\pi$ we have that
$\what{\til{\iota}\kappa}(\sig)=\tiljrpi\hat{\kappa}(\sig)$, by a computation
analagous to (\ref{eq:rpiisquot}), above.  Next, we wish to establish that
\begin{equation}\label{eq:gampiiisiigampih}
\gampi\comp\til{\iota}=(\til{\iota}_\pi\tens\til{\iota}_\pi)\comp\gampih.
\end{equation}
If $\kappa\in\meash$ and $x\in\bdophil\htens\bdophil$, then
\begin{align*}
\dpair{\gampi(\til{\iota}\kappa)}{x}
&=\int_G \dpair{\pi(s)\tens\pi(s)^*}{x}d\til{\iota}\kappa(s) \\
&=\int_H \dpair{\pi(s)\tens\pi(s)^*}{x}d\kappa(s)
=\dpair{\gampih(\kappa)}{x}
\end{align*}
whence, as elements of $\bdophil\ehtens\bdophil$, 
$\gampi(\til{\iota}\kappa)=\gampih(\kappa)$.  However, in $\bdophil$, the 
inclusion map $\mstarpih\hookrightarrow\mstarpi$ is the map $\til{\iota}_\pi$,
and thus (\ref{eq:gampiiisiigampih}) holds.

Now, since $\mu$ is supported on $H$, we have that $\mu=\til{\iota}\kappa$
for some $\kappa\in\meash$.  Then for each $(\sig,\tau)\iin E_\pi\cross E_\pi$
we have that
\begin{align*}
\hat{\mu}(\sig\tau^{-1})
&=\what{\til{\iota}\kappa}(\sig\tau^{-1})=\hat{\kappa}(r(\sig\tau^{-1})
=\hat{\kappa}(r_\pi(\sig)r_\pi(\tau^{-1})) \\
&=(\tilfpih\tens\tilfpih)\gampih(\kappa)(r_\pi(\sig),r_\pi(\tau)),\quad
\text{by part I} \\
&=(\tiljrpi\tens\tiljrpi)(\tilfpih\tens\tilfpih)\gampih(\kappa)(\sig,\tau) \\
&=(\tilfpi\tens\tilfpi)(\til{\iota}_\pi\tens\til{\iota}_\pi)
\gampih(\kappa)(\sig,\tau),\quad\text{by (\ref{eq:jfpihisfpii})} \\
&=(\tilfpi\tens\tilfpi)\gampi(\til{\iota}\kappa)(\sig,\tau),
\quad\text{by (\ref{eq:gampiiisiigampih})} \\
&=(\tilfpi\tens\tilfpi)\gampi(\mu)(\sig,\tau).
\end{align*}

{\bf III.}  We now cover the case of a general $\mu$ in $\measg$.

Let $U$ be a relatively compact symmetric open neighbourhood of
the identity in $G$.  Then $H=\bigcup_{n=1}^\infty U^n$ is a compactly
generated open subgroup of $G$.  We note that if $T$ is
a transversal for $H$ in $G$ then
\[
\mu=\sum_{t\in T}\mu_{tH}
\]
which is an absolutely summable series.  For each $t\iin T$ let 
\[
\mu_t=\del_{t^{-1}}\con(\mu_{tH})
\]
so $\supp{\mu_t}\subset H$ and $\mu=\sum_{t\in T}\del_t\con \mu_t$.
We then have that
\begin{align*}
(\til{F}_\pi\tens \til{F}_\pi)\gampi(\mu)
&=\sum_{t\in T}(\til{F}_\pi\tens \til{F}_\pi)\gampi(\del_t\con \mu_t) \\
&=\sum_{t\in T}(\til{F}_\pi\tens \til{F}_\pi)
          \left[(\pi(t)\tens\pi(t)^*)\gampi(\mu_t)\right] \\
&=\sum_{t\in T}(\hat{t}|_{E_\pi}\tens\bar{\hat{t}}|_{E_\pi})
(\til{F}_\pi\tens \til{F}_\pi)\gampi(\mu_t).
\end{align*}
Hence if $(\sig,\tau)\in E_\pi\cross E_\pi$, we obtain
\begin{align*}
(\til{F}_\pi\tens \til{F}_\pi)\gampi(\mu)(\sig,\tau)
&=\sum_{t\in T}\hat{t}(\sig)\wbar{\hat{t}(\tau)}
(\til{F}_\pi\tens \til{F}_\pi)\gampi(\mu_t)(\sig,\tau) \\
&=\sum_{t\in T}\hat{t}(\sig)\wbar{\hat{t}(\tau)}\hat{\mu}_t(\sig\tau^{-1}),
\qquad\text{by part II} \\
&=\sum_{t\in T}\hat{\del}_t(\sig\tau^{-1})\hat{\mu}_t(\sig\tau^{-1}) \\
&=\sum_{t\in T}\what{\del_t\con \mu_t}(\sig\tau^{-1})
=\hat{\mu}(\sig\tau^{-1}).
\end{align*}
Thus our first claim is established in general.

If $E$ is any closed subset of $\dgrp$ then by \cite[33.7]{hewittrII}
there is a representation $\pi:G\to\unophil$ for which $E_\pi=E$.  Hence
\[
u(\sig,\tau)=\hat{\mu}(\sig\tau^{-1})
=(\til{F}_\pi\tens \til{F}_\pi)\gampi(\mu)(\sig,\tau)
\]
defines an element of $\valb{E}$, and of $\valo{E}$ if $\mu\in\bloneg$.  
If $E$ is compact, we note that any representation $\pi$
for which $E_\pi=E$, is norm continuous by Proposition \ref{prop:repwunit},
since $\cstarpi\cong\conto{E}$, which is unital.  Thus $\gampi(\measg)=
\gampi(\bloneg)\subset\cstarpi\htens\cstarpi$.  Hence $u$, as above, 
is in $\valoc{E}$.  \endpf

We can now obtain a generalisation of \cite[Prop.\ 5.7]{stormer}.
This is a straightforward corollary of Theorems \ref{theo:posabelian}
and \ref{theo:abelianreps}.

\begin{abelianreps1}
If $G$ is a locally compact abelian group, $\pi:G\to\unophil$
is a strongly continuous representation and $\mu\in\measg$, then
the following are equivalent:

{\bf (i)} $\gampi(\mu)$ is positive.

{\bf (ii)} $\gampi(\mu)$ is completely positive.

{\bf (iii)} $(\sig,\tau)\mapsto\hat{\mu}(\sig\tau^{-1})$ is
positive definite on $E_\pi\cross E_\pi$.
\end{abelianreps1}

The next result follows directly from Theorem \ref{theo:abelianreps}, but
can also be deduced from Proposition \ref{prop:kergampi}.

\begin{abelianreps2}\label{cor:abelianreps2}
If $G$ is a locally compact abelian group and $\pi:G\to\unophil$
is a strongly continuous representation then
\[
\ker\gampi=\{\mu\in\measg:\hat{\mu}|_{E_\pi E_\pi^{-1}}=0\}.
\]
\end{abelianreps2}

Let us now address the assumption that $f\in\bloneg\setdif\augidealg$
in Theorem \ref{theo:repinhtp} (iii).  We want to show that having
an $f\iin\augidealg$ for which $\gampi(f)\subset\cstarpi\htens\cstarpi$
does not imply that $\pi$ is norm continuous.
First, following Corollary \ref{cor:abelianreps2} we see that
that $\gampi(f)=0$ if the support of $\hat{f}$
misses the difference set $E_\pi E_\pi^{-1}$.  Thus it is possible that
$\gampi(f)=0\in\cstarpi\htens\cstarpi$, though $\pi$ need not be
norm continuous, i.e.\ $E_\pi$ need not be compact.  Thus we may ask if
$\gampi(\bloneg)\cap\brac{\cstarpi\htens\cstarpi}=\{0\}$ when
$\pi$ is not norm continuous.  However, this may not happen, as the next 
example shows.

\begin{squareex}{\rm
Let $G=\Tee$, and identify $\what{\Tee}=\Zee$.
Define $\pi:\Tee\to\unop{\ltwo{\En}}$ for
each $z\iin\Tee$ by
\[
\pi(z)\bigl(\xi_n\bigr)_{n\in\En}=\brac{z^{n^2}\xi_n}_{n\in\En}.
\]
Then $E_\pi=\{n^2:n\in\En\}$, which is not compact in $\Zee$.
Hence $\cstarpi\cong\ceeo{E_\pi}$, which is not unital, so $\pi$ is not 
norm continuous on $\Tee$, by Proposition \ref{prop:repwunit}.  
Fix $k\iin\Zee\setdif\{0\}$ and let $\hat{k}(z)=z^k$. 
Then for each pair $n,\;m$ in $\En$, 
using normalized Haar measure on $\Tee$ and Theorem 
\ref{theo:abelianreps}, we have that
\[
(F_\pi\tens F_\pi)\gampi(\hat{k})(\hat{n},\hat{m})
=\int_\Tee z^kz^{n^2}\bar{z}^{m^2}dz
=\begin{cases} 1 & \iif m^2-n^2=k \\ 0 &\text{ otherwise}  \end{cases}.
\]
The set of solutions to $m^2-n^2=(m-n)(m+n)=k$ is clearly finite;
we shall write them $\{(n_1,m_1),\dots,(n_{l(k)},m_{l(k)})\}$.
We then see that
\[
(F_\pi\tens F_\pi)\gampi(\hat{k})=\sum_{i=1}^{l(k)}1_{(n_i,m_i)}
=\sum_{i=1}^{l(k)}1_{n_i}\tens 1_{m_i}\in\valoc{E_\pi}.
\]
Hence $\gampi(\hat{k})\in\cstarpi\htens\cstarpi$.
In fact, since $\augideal{\Tee}=\wbar{\spn}
\cbrac{\hat{k}:k\in\Zee\setdif\{0\}}$, we have that
$\gampi(\augideal{\Tee})\subset\cstarpi\htens\cstarpi$.}
\end{squareex}

We remark that for a general locally compact abelian group $G$, 
and representation $\pi:G\to\unophil$, $\valoc{E_\pi}\subset 
\conto{E_\pi\cross E_\pi}$.  Thus
if $f\iin\bloneg$ is such that
$\hat{f}(\sig)\not=0$ for some $\sig\iin\dgrp$ such that $E_\pi\cap\sig E_\pi$
is not compact, then $\gampi(f)\not\in\cstarpi\htens\cstarpi$ by Theorem
\ref{theo:abelianreps}.  Thus if $E_\pi\cap\sig E_\pi$ is compact for no
$\sig\iin\dgrp$ then we have that
\[
\gampi(\bloneg)\cap\brac{\cstarpi\htens\cstarpi}=\{0\}.
\]
Note that $E_\pi\cap\sig E_\pi$ is never compact if $E_\pi=\dgrp$, 
which occurs, for example when $\pi$ is the left regular representation 
$\lam$. It is shown in \cite[Cor.\ 4.7]{stormer} 
that $\Gam_\lam$ is an isometry.  This was extended to non-abelian groups
in  \cite{ghahramani} and expanded upon in \cite{neufang}, while 
\cite{neufangs} contains a proof that $\Gam_\lam$ is a complete isometry.
An analogue for the Fourier algebra of an amenable group
is shown in \cite[Cor.\ 5.4]{spronkt}. 

\begin{injectiveex}{\rm
If $G=\Ree$, then $[0,\infty)\cap\bigl(
s+[0,\infty)\bigr)=[\min\{s,0\},\infty)$ is never compact.  Thus if
$\pi$ is a representation of $\Ree$ such that $E_\pi=[0,\infty)$, then
$\gampi:\meas{\Ree}\to\mstarpi\ehtens\mstarpi$ is injective by 
Corollary \ref{cor:abelianreps2}.  Is $\gampi$ isometric? 
How about $\gampi|_{\blone{\Ree}}$?
More generally, under what conditions for an arbitrary abelian group $G$ and
representation $\pi$ is $\gampi$, or $\gampi|_{\bloneg}$, a quotient map?}
\end{injectiveex}

\nocite{sakai}

{
\bibliography{normcontbib}
\bibliographystyle{plain}
}

\smallskip
{\sc Department of Mathematics, Texas A\& M University, College Station, 
Texas 77843-3368, U.S.A.}

{\it E-mail address:} {\tt rsmith@math.tamu.edu, spronk@math.tamu.edu}

\end{document}